\documentclass[preprint,12pt]{elsarticle}
\usepackage{amssymb}
\usepackage{amsmath}

\newtheorem{thm}{Theorem}[section]
\newtheorem{lem}[thm]{Lemma}
\newtheorem{cor}[thm]{Corollary}
\newdefinition{defn}[thm]{Definition}
\newdefinition{rem}[thm]{Remark}
\newproof{pf}{Proof}
\newproof{pot1}{Proof of Theorem \ref{czx}}
\newproof{pot2}{Proof of Theorem \ref{wqj}}

\journal{arXiv}
\begin{document}
\begin{frontmatter}

\title{Infinitely Many Weak Solutions for Fractional Dirichlet Problem with $p$-Laplacian}

\author{Taiyong Chen,\ \ Wenbin Liu
\footnote{Corresponding author.\\
{\it Telephone number:} (86-516) 83591530. {\it Fax number:} (86-516) 83591591.
{\it E-mail addresses:} taiyongchen@cumt.edu.cn (T. Chen), wblium@163.com (W. Liu), jinhua197927@163.com (H. Jin).}
,\ \ Hua Jin}

\address{Department of Mathematics, China University of Mining and Technology, Xuzhou 221116, PR China}

\begin{abstract}
We focus on the study of $p$-Laplacian Dirichlet problem containing the left and right fractional derivative operators. By using the genus properties in critical point theory, we establish some new criteria to guarantee the existence of infinitely many weak solutions for the considered problem.
\end{abstract}


\end{frontmatter}

\section{Introduction}
\label{sec1}
The purpose of the present paper is to study the existence of infinitely many weak solutions for the fractional Dirichlet problem with $p$-Laplacian of the form
\begin{eqnarray}
\label{dbvp}
\left\{
\begin{array}{ll}
_tD_T^\alpha\phi_p({_0D_t^\alpha}u(t))=f(t,u(t)),\ \ t\in(0,T),\\
u(0)=u(T)=0,
\end{array}
\right.
\end{eqnarray}
where $p>1$ is a constant, $_0D_t^\alpha$ and $_tD_T^\alpha$ are the left and right Riemann-Liouville fractional derivatives of order $\alpha\in(1/p,1]$ respectively, $\phi_p:\mathbb{R}\rightarrow\mathbb{R}$ is the $p$-Laplacian defined by
\begin{eqnarray*}
\phi_p(s)=|s|^{p-2}s\ (s\neq0),\ \ \phi_p(0)=0,
\end{eqnarray*}
and $f\in C([0,T]\times\mathbb{R},\mathbb{R})$.

Since the fractional derivatives provide an excellent tool to describe the memory and hereditary properties of various materials and processes (see \cite{dkfa,3,hilf,kjfx,main}), the fractional order models are more adequate than the integer order models in some real world problems. Moreover the $p$-Laplacian introduced by Leibenson (see \cite{lsl}) often occurs in non-Newtonian fluid theory, nonlinear elastic mechanics and so forth. Note that, when $p=2$, the nonlinear and nonlocal differential operator $_tD_T^\alpha\phi_p({_0D_t^\alpha})$ reduces to the linear differential operator $_tD_T^\alpha{_0D_t^\alpha}$, and further reduces to the local second-order differential operator $-d^2/dt^2$ when $\alpha=1$.

In the past decade, many results on the existence and multiplicity of solutions for the nonlinear fractional boundary value problems (BVPs for short) have been obtained (see \cite{21,zbhl,23,e1,25,wj}). In addition, as the applications in physical phenomena exhibiting anomalous diffusion, the models containing left and right fractional differential operators are recently gaining more attention. We refer the readers to \cite{jbe3,jcre,jerw,12,fjy,zzr}.

In \cite{bvp3} the authors considered the existence of solutions for BVP (\ref{dbvp}) when $F(t,x)$ satisfies the so-called Ambrosetti-Rabinowtiz condition. More explicitly, the nonlinearity $f$ is supposed to satisfy the following assumption.

$(H)$ There exist constants $\mu\in(0,1/p)$ and $M>0$ such that
\begin{eqnarray*}
0<F(t,x)\leq\mu xf(t,x),\ \ \forall t\in[0,T],\ x\in\mathbb{R}\ \mbox{with}\ |x|\geq M,
\end{eqnarray*}
where $F(t,x)=\int_0^xf(t,s)ds$. \\
Obviously the assumption $(H)$ implies
\begin{eqnarray*}
F(t,x)\geq c_1|x|^{\frac{1}{\mu}}-c_2,\ \ \forall (t,x)\in[0,T]\times\mathbb{R},
\end{eqnarray*}
where $c_1,c_2>0$ are two constants.

Motivated by the above works, the aim of this paper is to study the existence of infinitely many weak solutions for BVP (\ref{dbvp}) under some assumptions different from $(H)$. The main ingredient used here is the genus properties in critical point theory.

In order to state our main results, we make some assumptions on the nonlinearity $f$ as follows.

$(H_1)$ There exist a constant $1<r_1<p$ and a function $a\in L^1([0,T],\mathbb{R}^+)$ such that
\begin{eqnarray*}
|f(t,x)|\leq r_1a(t)|x|^{r_1-1},\ \ \forall (t,x)\in[0,T]\times\mathbb{R}.
\end{eqnarray*}

$(H_2)$ There exist an open interval $\mathbb{I}\subset[0,T]$ and three constants $\eta,\delta>0$, $1<r_2<p$ such that
\begin{eqnarray*}
F(t,x)\geq\eta|x|^{r_2},\ \ \forall (t,x)\in\mathbb{I}\times[-\delta,\delta].
\end{eqnarray*}

$(H_3)$ $f(t,x)=-f(t,-x),\ \ \forall (t,x)\in[0,T]\times\mathbb{R}$.

We are now in a position to state our main results.

\begin{thm}
\label{czx}
Let $(H_1)$ and $(H_2)$ be satisfied. Then BVP (\ref{dbvp}) possesses at least one nontrivial weak solution.
\end{thm}

\begin{thm}
\label{wqj}
Let $(H_1)$-$(H_3)$ be satisfied. Then BVP (\ref{dbvp}) possesses infinitely many nontrivial weak solutions.
\end{thm}

\begin{rem}
It is easy to verify from $(H_1)$ that
\begin{eqnarray*}
|F(t,x)|\leq a(t)|x|^{r_1},\ \ \forall (t,x)\in[0,T]\times\mathbb{R},
\end{eqnarray*}
which is obviously different from the assumption $(H)$. In addition, in $(H)$, $F(t,x)$ is supposed to be nonnegative when $|x|$ is big enough. However $F(t,x)$ in the present paper can change its sign.
\end{rem}

The rest of this paper is organized as follows. Section \ref{sec2} contains some preliminary results. In Section \ref{sec4}, the proof of main results is given.

\section{Preliminaries}
\label{sec2}

\subsection{Fractional Sobolev space}
In this subsection, some definitions and notations of the fractional calculus are presented (see \cite{15,18}). Moreover we introduce a fractional Sobolev space and some properties of this space (see \cite{fjy}).

\begin{defn}%[\cite{15}]
\label{defn2.1}
For $\gamma>0$, the left and right Riemann-Liouville fractional integrals of order $\gamma$ of a function $u:[a,b]\rightarrow\mathbb{R}$ are given by
\begin{eqnarray*}
&&_aI_{t}^\gamma u(t)=\frac{1}{\Gamma(\gamma)}\int_a^t(t-s)^{\gamma -1}u(s)ds,\\
&&_tI_{b}^\gamma u(t)=\frac{1}{\Gamma(\gamma)}\int_t^b(s-t)^{\gamma -1}u(s)ds,
\end{eqnarray*}
provided that the right-hand side integrals are pointwise defined on $[a,b]$, where $\Gamma(\cdot)$ is the Gamma function.
\end{defn}

\begin{defn}%[\cite{15}]
\label{defn2.2}
For $n-1\leq\gamma<n\ (n\in\mathbb{N})$, the left and right Riemann-Liouville fractional derivatives of order $\gamma$ of a function $u:[a,b]\rightarrow\mathbb{R}$ are given by
\begin{eqnarray*}
&&_aD_{t}^\gamma u(t)=\frac{d^n}{dt^n}{_a}I_{t}^{n-\gamma} u(t),\\
&&_tD_{b}^\gamma u(t)=(-1)^n\frac{d^n}{dt^n}{_t}I_{b}^{n-\gamma} u(t).
\end{eqnarray*}
\end{defn}

\begin{rem}
\label{zxzj}
When $\gamma=1$, we can obtain from Definition \ref{defn2.1} and \ref{defn2.2} that
\begin{eqnarray*}
_aD_{t}^1 u(t)=u'(t),\ \ _tD_{b}^1 u(t)=-u'(t),
\end{eqnarray*}
where $u'$ is the usual first-order derivative of $u$.
\end{rem}

\begin{defn}%[\cite{cjh}]
\label{defn3.1}
For $0<\alpha\leq1$ and $1<p<\infty$, the fractional derivative space $E{_0^{\alpha,p}}$ is defined by the closure of $C_0^\infty([0,T],\mathbb{R})$ with respect to the following norm
\begin{eqnarray*}
\|u\|_{E^{\alpha,p}}=(\|u\|_{L^p}^p+\|{_0D_t^\alpha}u\|_{L^p}^p)^{\frac{1}{p}},\ \ \forall u\in E{_0^{\alpha,p}},
\end{eqnarray*}
where $\|u\|_{L^p}=\left(\int_0^T|u(t)|^pdt\right)^{1/p}$ is the norm of $L^p([0,T],\mathbb{R})$.
\end{defn}

\begin{rem}
It is obvious that, for $u\in E{_0^{\alpha,p}}$, one has
\begin{eqnarray*}
u,{_0D_t^\alpha}u\in L^p([0,T],\mathbb{R}),\ \ u(0)=u(T)=0.
\end{eqnarray*}
\end{rem}

\begin{lem}[see \cite{fjy}]
\label{lem1}
Let $0<\alpha\leq1$ and $1<p<\infty$. The fractional derivative space $E{_0^{\alpha,p}}$ is a reflexive and separable Banach space.
\end{lem}

\begin{lem}[see \cite{fjy}]
Let $0<\alpha\leq1$ and $1<p<\infty$. For $u\in E_0^{\alpha,p}$, one has
\begin{eqnarray}
\label{cp}
\|u\|_{L^p}\leq C_p\|{_0D_t^\alpha}u\|_{L^p},
\end{eqnarray}
where
\begin{eqnarray*}
C_p=\frac{T^\alpha}{\Gamma(\alpha+1)}>0
\end{eqnarray*}
is a constant. Moreover, if $\alpha>1/p$, then
\begin{eqnarray}
\label{cwq}
\|u\|_\infty\leq C_\infty\|u\|_{E^{\alpha,p}},
\end{eqnarray}
where $\|u\|_\infty=\max_{t\in[0,T]}|u(t)|$ is the norm of $C([0,T],\mathbb{R})$ and
\begin{eqnarray*}
C_\infty=\frac{T^{\alpha-\frac{1}p}}{\Gamma(\alpha)(\alpha q-q+1)^{\frac{1}q}}>0,\ \
q=\frac{p}{p-1}>1
\end{eqnarray*}
are two constants.
\end{lem}

\begin{rem}
From (\ref{cp}), we know that the norm $\|\cdot\|_{E^{\alpha,p}}$ of $E_0^{\alpha,p}$ is equivalent to the norm $\|{_0D_t^\alpha}\cdot\|_{L^p}$. Hence we can consider the space $E_0^{\alpha,p}$ with norm $\|{_0D_t^\alpha}\cdot\|_{L^p}$ in the following analysis.
\end{rem}

\begin{lem}[see \cite{fjy}]
\label{thm3.2}
Let $1/p<\alpha\leq1$ and $1<p<\infty$. The imbedding of $E_0^{\alpha,p}$ in $C([0,T],\mathbb{R})$ is compact.
\end{lem}

\subsection{Critical point theory}

Now we introduce some notations and necessary definitions of the critical point theory (see \cite{dl1,dl2}). Let $X$ be a real Banach space, $I\in C^1(X,\mathbb{R})$ which means that $I$ is a continuously Fr\'{e}chet differentiable functional defined on $X$.

\begin{defn}%[see \cite{dl1}]
Let $I\in C^1(X,\mathbb{R})$. If any sequence $\{u_k\}\subset X$ for which $\{I(u_k)\}$ is bounded and $I'(u_k)\rightarrow0$ as $k\rightarrow\infty$ possesses a convergent subsequence in $X$, then we say that $I$ satisfies the Palais-Smale condition ((PS)-condition for short).
\end{defn}

\begin{lem}[see \cite{dl1}]
\label{dll1}
Let $X$ be a real Banach space and $I\in C^1(X,\mathbb{R})$ satisfies the (PS)-condition. If $I$ is bounded from blow, then $c=\inf_X I$ is a critical value of $I$.
\end{lem}

In order to find infinitely many nontrivial critical points of $I$, the following {\it genus} properties are needed in our argument. Set
\begin{eqnarray*}
&&\Sigma=\{A\subset X-\{0\}|A\mbox{ is closed in }X\mbox{ and symmetric with respect to }0\},\\
&&K_c=\{u\in X|I(u)=c,\ I'(u)=0\},\ \ I^c=\{u\in X|I(u)\leq c\}.
\end{eqnarray*}

\begin{defn}%[see \cite{dl2}]
For $A\in\Sigma$, we say the genus of $A$ is $n$ denoted by $\gamma(A)=n$ if there is an odd map $G\in C(A,\mathbb{R}^n\backslash\{0\})$ and $n$ is the smallest integer with this property.
\end{defn}

\begin{lem}[see \cite{dl2}]
\label{dll2}
Let $I$ be an even $C^1$ functional on $X$ and satisfy the (PS)-condition. For any $n\in\mathbb{N}$, set
\begin{eqnarray*}
\Sigma_n=\{A\in\Sigma|\gamma(A)\geq n\},\ \ c_n=\inf_{A\in\Sigma_n}\sup_{u\in A}I(u).
\end{eqnarray*}
(i) If $\Sigma_n\neq\emptyset$ and $c_n\in\mathbb{R}$, then $c_n$ is a critical value of $I$.\\
(ii) If there exists $l\in\mathbb{N}$ such that $c_n=c_{n+1}=\cdots=c_{n+l}=c\in\mathbb{R}$, and $c\neq I(0)$, then $\gamma(K_c)\geq l+1$.
\end{lem}

\begin{rem}
From Remark 7.3 in \cite{dl2}, we know that if $K_c\in\Sigma$ and $\gamma(K_c)>1$, then $K_c$ contains infinitely many distinct points, that is, $I$ has infinitely many distinct critical points in $X$.
\end{rem}

\section{Proof of main results}
\label{sec4}

The aim of this section is to prove our main results. For this purpose, we are going to set up the corresponding variational framework to obtain the solutions of BVP (\ref{dbvp}).

Define the functional $I:E_0^{\alpha,p}\rightarrow\mathbb{R}$ by
\begin{align}
\label{fh}
I(u)
&=\int_0^T\left(\frac{1}{p}|{_0D_t^\alpha}u(t)|^p-F(t,u(t))\right)dt\nonumber\\
&=\frac{1}{p}\|u\|_{E^{\alpha,p}}^p-\int_0^TF(t,u(t))dt.
\end{align}
Throughout this paper, by the weak solutions of BVP (\ref{dbvp}) we mean the critical points of the associated energy functional $I$. It is easy to verify from $f\in C([0,T]\times\mathbb{R},\mathbb{R})$ and (\ref{cwq}) that the functional $I$ is well defined on $E_0^{\alpha,p}$ and is a continuously Fr\'{e}chet differentiable functional, that is, $I\in C^1(E_0^{\alpha,p},\mathbb{R})$. Moreover we have
\begin{eqnarray}
\label{fhds}
\langle I'(u),v\rangle
=\int_0^T\phi_p({_0D_t^\alpha}u(t)){_0D_t^\alpha}v(t)dt
-\int_0^Tf(t,u(t))v(t)dt
\end{eqnarray}
for all $u,v\in E_0^{\alpha,p}$.

\begin{lem}
\label{xyj}
Suppose that $(H_1)$ is satisfied. Then $I$ is bounded from below in $E_0^{\alpha,p}$.
\end{lem}

\begin{pf}
From $(H_1)$, one has
\begin{eqnarray*}
\label{xtj}
|F(t,u)|\leq a(t)|u|^{r_1},\ \ \forall (t,u)\in[0,T]\times\mathbb{R},
\end{eqnarray*}
which together with (\ref{cwq}) and (\ref{fh}) yields
\begin{align}
\label{fsdy}
I(u)
&\geq\frac{1}{p}\|u\|_{E^{\alpha,p}}^p-\int_0^Ta(t)|u(t)|^{r_1}dt\nonumber\\
&\geq\frac{1}{p}\|u\|_{E^{\alpha,p}}^p-\|a\|_{L^1}\|u\|_\infty^{r_1}\nonumber\\
&\geq\frac{1}{p}\|u\|_{E^{\alpha,p}}^p-C_\infty^{r_1}\|a\|_{L^1}\|u\|_{E^{\alpha,p}}^{r_1}.
\end{align}
Since $1<r_1<p$, (\ref{fsdy}) implies $I(u)\rightarrow\infty$ as $\|u\|_{E^{\alpha,p}}\rightarrow\infty$. So $I$ is bounded from below. $\Box$
\end{pf}

\begin{lem}
\label{ps}
Assume that $(H_1)$ holds. Then $I$ satisfies the (PS)-condition in $E_0^{\alpha,p}$.
\end{lem}

\begin{pf}
Let $\{u_k\}\subset E_0^{\alpha,p}$ be a sequence such that
\begin{eqnarray*}
|I(u_k)|\leq K,\ \ I'(u_k)\rightarrow0\ \mbox{as}\ {k\rightarrow\infty},
\end{eqnarray*}
where $K>0$ is a constant. Then (\ref{fsdy}) implies that $\{u_k\}$ is bounded in $E_0^{\alpha,p}$. Since $E_0^{\alpha,p}$ is a reflexive Banach space (see Lemma \ref{lem1}), going if necessary to a subsequence, we can assume $u_k\rightharpoonup u$ in $E_0^{\alpha,p}$. Hence we obtain from $I'(u_k)\rightarrow0\ \mbox{as}\ {k\rightarrow\infty}$ and the definition of weak convergence that
\begin{align}
\label{jl3}
&\langle I'(u_k)-I'(u),u_k-u\rangle \nonumber\\
&=\langle I'(u_k),u_k-u\rangle -\langle I'(u),u_k-u\rangle \nonumber\\
&\leq\|I'(u_k)\|_{(E_0^{\alpha,p})^*}\|u_k-u\|_{E^{\alpha,p}}-\langle I'(u),u_k-u\rangle \nonumber\\
&\rightarrow0\ \ \mbox{as}\ k\rightarrow\infty.
\end{align}
On the other hand, one can derive from (\ref{cwq}) and Lemma \ref{thm3.2} that $\{u_k\}$ is bounded in $C([0,T],\mathbb{R})$ and $\|u_k-u\|_\infty\rightarrow0$ as $k\rightarrow\infty$. So there exists a constant $C_1>0$ such that
\begin{eqnarray*}
|f(t,u_k(t))-f(t,u(t))|\leq C_1,\ \ \forall t\in[0,T],
\end{eqnarray*}
which yields
\begin{align}
\label{jl4}
\left|\int_0^T(f(t,u_k(t))-f(t,u(t)))(u_k(t)-u(t))dt\right|
&\leq C_1T\|u_k-u\|_\infty\nonumber\\
&\rightarrow0\ \ \mbox{as}\ k\rightarrow\infty.
\end{align}
By (\ref{fhds}), one gets
\begin{align*}
&\langle I'(u_k)-I'(u),u_k-u\rangle \\
&=\int_0^T(\phi_p({_0D_t^\alpha}u_k(t))-\phi_p({_0D_t^\alpha}u(t)))
({_0D_t^\alpha}u_k(t)-{_0D_t^\alpha}u(t))dt\\
&\ \ \ \ -\int_0^T(f(t,u_k(t))-f(t,u(t)))(u_k(t)-u(t))dt,
\end{align*}
which together with (\ref{jl3}) and (\ref{jl4}) yields
\begin{eqnarray}
\label{jl8}
\int_0^T(\phi_p({_0D_t^\alpha}u_k(t))-\phi_p({_0D_t^\alpha}u(t)))
({_0D_t^\alpha}u_k(t)-{_0D_t^\alpha}u(t))dt
\rightarrow0
\end{eqnarray}
as $k\rightarrow\infty$.

Following (2.10) in \cite{js}, we can find the constants $C_2,C_3>0$ such that
\begin{align}
\label{jl5}
&\int_0^T(\phi_p({_0D_t^\alpha}u_k(t))-\phi_p({_0D_t^\alpha}u(t)))
({_0D_t^\alpha}u_k(t)-{_0D_t^\alpha}u(t))dt\nonumber\\
&\geq
\left\{
\begin{array}{ll}
C_2\int_0^T|{_0D_t^\alpha}u_k(t)-{_0D_t^\alpha}u(t)|^pdt,\ \ p\geq2,\\
C_3\int_0^T\frac{|{_0D_t^\alpha}u_k(t)-{_0D_t^\alpha}u(t)|^2}
{(|{_0D_t^\alpha}u_k(t)|+|{_0D_t^\alpha}u(t)|)^{2-p}}dt,\ \ 1<p<2.
\end{array}
\right.
\end{align}
When $1<p<2$, by the H\"older inequality, one has
\begin{align*}
&\int_0^T|{_0D_t^\alpha}u_k(t)-{_0D_t^\alpha}u(t)|^pdt\\
&\leq\left(\int_0^T\frac{|{_0D_t^\alpha}u_k(t)-{_0D_t^\alpha}u(t)|^2}
{(|{_0D_t^\alpha}u_k(t)|+|{_0D_t^\alpha}u(t)|)^{2-p}}dt\right)^{\frac{p}{2}}\\
&\ \ \ \ \cdot\left(\int_0^T(|{_0D_t^\alpha}u_k(t)|+|{_0D_t^\alpha}u(t)|)^pdt\right)
^{\frac{2-p}{2}}\\
&\leq C_4(\|u_k\|_{E^{\alpha,p}}^p+\|u\|_{E^{\alpha,p}}^p)^{\frac{2-p}{2}}
\left(\int_0^T\frac{|{_0D_t^\alpha}u_k(t)-{_0D_t^\alpha}u(t)|^2}
{(|{_0D_t^\alpha}u_k(t)|+|{_0D_t^\alpha}u(t)|)^{2-p}}dt\right)^{\frac{p}{2}},
\end{align*}
where $C_4=2^{(p-1)(2-p)/2}>0$ is a constant, which together with (\ref{jl5}) implies
\begin{align}
\label{jl6}
&\int_0^T(\phi_p({_0D_t^\alpha}u_k(t))-\phi_p({_0D_t^\alpha}u(t)))
({_0D_t^\alpha}u_k(t)-{_0D_t^\alpha}u(t))dt\nonumber\\
&\geq C_3C_4^{-\frac{2}{p}}(\|u_k\|_{E^{\alpha,p}}^p+\|u\|_{E^{\alpha,p}}^p)^{\frac{p-2}{p}}
\|u_k-u\|_{E^{\alpha,p}}^2,\ \ 1<p<2.
\end{align}
When $p\geq2$, by (\ref{jl5}), we get
\begin{align}
\label{jl7}
&\int_0^T(\phi_p({_0D_t^\alpha}u_k(t))-\phi_p({_0D_t^\alpha}u(t)))
({_0D_t^\alpha}u_k(t)-{_0D_t^\alpha}u(t))dt\nonumber\\
&\geq C_2\|u_k-u\|_{E^{\alpha,p}}^p,\ \ p\geq2.
\end{align}
Then it follows from (\ref{jl8}), (\ref{jl6}) and (\ref{jl7}) that
\begin{eqnarray*}
\|u_k-u\|_{E^{\alpha,p}}\rightarrow0\ \ \mbox{as}\ \ k\rightarrow\infty,
\end{eqnarray*}
which means that $I$ satisfies the (PS)-condition. $\Box$
\end{pf}

\begin{pot1}
From Lemma \ref{dll1}, \ref{xyj} and \ref{ps}, we know that $c=\inf_{E_0^{\alpha,p}}I(u)$ is a critical value of $I$, that is, there exists a critical point $u^*\in E_0^{\alpha,p}$ such that $I(u^*)=c$.

Finally we prove $u^*\neq0$. Let $u_0\in\left(W_0^{1,2}(\mathbb{I},\mathbb{R})\cap E_0^{\alpha,p}\right)\setminus\{0\}$ and $\|u_0\|_\infty=1$, by (\ref{fh}) and $(H_2)$, we have
\begin{align}
\label{jfl}
I(su_0)
&=\frac{1}{p}\|su_0\|_{E^{\alpha,p}}^p-\int_0^TF(t,su_0(t))dt\nonumber\\
&=\frac{s^p}{p}\|u_0\|_{E^{\alpha,p}}^p-\int_\mathbb{I}F(t,su_0(t))dt\nonumber\\
&\leq\frac{s^p}{p}\|u_0\|_{E^{\alpha,p}}^p-\eta s^{r_2}\int_\mathbb{I}|u_0(t)|^{r_2}dt,\ \ 0<s\leq \delta.
\end{align}
Since $1<r_2<p$, it follows from (\ref{jfl}) that $I(su_0)<0$ for $s>0$ small enough. Thus $I(u^*)=c<0$, hence $u^*$ is a nontrivial critical point of $I$, and so $u^*$ is a nontrivial solution of BVP (\ref{dbvp}). $\Box$
\end{pot1}

\begin{pot2}
From Lemma \ref{xyj} and \ref{ps}, we know that $I\in C^1(E_0^{\alpha,p},\mathbb{R})$ is bounded from below and satisfies the (PS)-condition. Moreover (\ref{fh}) and $(H_3)$ show that $I$ is even and $I(0)=0$.

Fixing $n\in\mathbb{N}$, we take $n$ disjoint open intervals $\mathbb{I}_i$ such that $\cup_{i=1}^n\mathbb{I}_i\subset\mathbb{I}$. Let $u_i\in\left(W_0^{1,2}(\mathbb{I}_i,\mathbb{R})\cap E_0^{\alpha,p}\right)\setminus\{0\}$ and $\|u_i\|_{E^{\alpha,p}}=1$, and
\begin{eqnarray*}
E_n=\mbox{span}\{u_1,u_2,\cdots,u_n\},\ \ S_n=\{u\in E_n|\|u\|_{E^{\alpha,p}}=1\}.
\end{eqnarray*}
For $u\in E_n$, there exist $\lambda_i\in\mathbb{R}$, such that
\begin{eqnarray}
\label{zks}
u(t)=\sum_{i=1}^n\lambda_iu_i(t),\ \ \forall t\in[0,T].
\end{eqnarray}
Then we have
\begin{align}
\label{nyfs}
\|u\|_{E^{\alpha,p}}^p
&=\int_0^T|_0D_t^\alpha u(t)|^pdt
=\sum_{i=1}^n|\lambda_i|^p\int_{\mathbb{I}_i}|_0D_t^\alpha u_i(t)|^pdt\nonumber\\
&=\sum_{i=1}^n|\lambda_i|^p\int_0^T|_0D_t^\alpha u_i(t)|^pdt
=\sum_{i=1}^n|\lambda_i|^p\|u_i\|_{E^{\alpha,p}}^p\nonumber\\
&=\sum_{i=1}^n|\lambda_i|^p,\ \ \forall u\in E_n.
\end{align}
By (\ref{cwq}), (\ref{fh}), (\ref{zks}) and $(H_2)$, for $u\in S_n$, one has
\begin{align}
\label{fnl}
I(su)
&=\frac{1}{p}\|su\|_{E^{\alpha,p}}^p-\int_0^TF(t,su(t))dt\nonumber\\
&=\frac{s^p}{p}\|u\|_{E^{\alpha,p}}^p-\sum_{i=1}^n\int_{\mathbb{I}_i}F(t,s\lambda_iu_i(t))dt\nonumber\\
&\leq\frac{s^p}{p}-\eta s^{r_2}\sum_{i=1}^n|\lambda_i|^{r_2}\int_{\mathbb{I}_i}|u_i(t)|^{r_2}dt
,\ \ 0<s\leq\frac{\delta}{C_\infty\lambda^*},
\end{align}
where $\lambda^*=\max_{i\in\{1,2,\cdots,n\}}|\lambda_i|>0$ is a constant. Since $1<r_2<p$, (\ref{fnl}) implies that there exist $\epsilon,\sigma>0$ such that
\begin{eqnarray}
\label{nlz}
I(\sigma u)<-\epsilon,\ \ \forall u\in S_n.
\end{eqnarray}
Let
\begin{eqnarray*}
S_n^\sigma=\{\sigma u|u\in S_n\},\ \ \Lambda=\left\{(\lambda_1,\lambda_2,\cdots,\lambda_n)\in\mathbb{R}^n|
\sum_{i=1}^n|\lambda_i|^p<\sigma^p\right\}.
\end{eqnarray*}
Thus we obtain from (\ref{nlz}) that
\begin{eqnarray*}
I(u)<-\epsilon,\ \ \forall u\in S_n^\sigma,
\end{eqnarray*}
which, together with the fact that $I$ is even and $I(0)=0$, yields
\begin{eqnarray*}
S_n^\sigma\subset I^{-\epsilon}\in\Sigma.
\end{eqnarray*}
On the other hand, from (\ref{nyfs}), we know that the mapping $(\lambda_1,\lambda_2,\cdots,\lambda_n)\rightarrow\sum_{i=1}^n\lambda_iu_i(t)$ from $\partial\Lambda$ to $S_n^\sigma$ is odd and homeomorphic. Hence, by some properties of the genus (see Proposition 7.5 and 7.7 in \cite{dl2}), one can deduce that
\begin{eqnarray*}
\gamma(I^{-\epsilon})\geq\gamma(S_n^\sigma)=n.
\end{eqnarray*}
Then $I^{-\epsilon}\in \Sigma_n$ and so $\Sigma_n\neq\emptyset$. Let
\begin{eqnarray*}
c_n=\inf_{A\in\Sigma_n}\sup_{u\in A}I(u).
\end{eqnarray*}
It follows from $I$ is bounded from below that $-\infty<c_n\leq-\epsilon<0$. That is, for any $n\in\mathbb{N}$, $c_n$ is a real negative number. Hence, by Lemma \ref{dll2}, $I$ admits infinitely many nontrivial critical points, and so BVP (\ref{dbvp}) possesses infinitely many nontrivial negative energy solutions. $\Box$
\end{pot2}

Obviously the following assumption $(H_2')$ implies $(H_2)$ holds.

$(H_2')$ There exist an open interval $\mathbb{I}\subset[0,T]$ and three constants $\eta,\delta>0$, $1<r<p$ such that
\begin{eqnarray*}
f(t,x)x\geq r\eta|x|^r,\ \ \forall (t,x)\in\mathbb{I}\times[-\delta,\delta].
\end{eqnarray*}
Moreover the assumption $(H_4)$ as follows implies $(H_1)$-$(H_3)$ hold.

$(H_4)$ $f(t,x)=rb(t)|x|^{r-2}x$, where $1<r<p$ is a constant, $b\in C([0,T],\mathbb{R})$ and there exists an open interval $\mathbb{I}\subset[0,T]$ such that $b(t)>0,\ \forall t\in\mathbb{I}$.\\
Then two natural corollaries can be stated as follows.

\begin{cor}
Let $(H_1)$ and $(H_2')$ be satisfied. Then BVP (\ref{dbvp}) possesses at least one nontrivial weak solution. Furthermore, if $(H_3)$ is also satisfied, then BVP (\ref{dbvp}) possesses infinitely many nontrivial weak solutions.
\end{cor}

\begin{cor}
Let $(H_4)$ be satisfied. Then BVP (\ref{dbvp}) possesses infinitely many nontrivial weak solutions.
\end{cor}

%\section*{Competing interests}
%The authors declare that they have no competing interests.

%\section*{Authors' contributions}
%The authors contributed equally in this article. They read and approved the final manuscript.

\section*{Acknowledgements}
This work was supported by the National Natural Science Foundation of China (11271364) and the Nature Science Foundation of Jiangsu Province (BK20130170).


\section*{References}
\begin{thebibliography}{0000}

\bibitem{21} R.P. Agarwal, D. O'Regan, S. Stanek,
Positive solutions for Dirichlet problems of singular nonlinear fractional differential equations,
{\it J. Math. Anal. Appl.} 371 (2010) 57-68.

\bibitem{zbhl} Z. Bai, H. L\"{u},
Positive solutions for boundary value problem of nonlinear fractional differential equation,
{\it J. Math. Anal. Appl.} 311 (2005) 495-505.

\bibitem{23} M. Benchohra, S. Hamani, S.K. Ntouyas,
Boundary value problems for differential equations with fractional order and nonlocal conditions,
{\it Nonlinear Anal.} 71 (2009) 2391-2396.

\bibitem{jbe3} D.A. Benson, S.W. Wheatcraft, M.M. Meerschaert,
The fractional-order governing equation of L\'{e}vy motion,
{\it Water Resour. Res.} 36 (2000) 1413-1423.

\bibitem{e1} G.M. Bisci, D. Repovs,
Higher nonlocal problems with bounded potential,
{\it J. Math. Anal. Appl.} 420 (2014) 167-176.

\bibitem{bvp3} T. Chen, W. Liu,
Solvability of fractional boundary value problem with $p$-Laplacian via critical point theory,
{\it Bound. Value Probl.} 2016 (2016) Article ID 75.

\bibitem{jcre} J. Cresson,
Inverse problem of fractional calculus of variations for partial differential equations,
{\it Commun. Nonlin. Sci. Numer. Simul.} 15 (2010) 987-996.

\bibitem{25} M.A. Darwish, S.K. Ntouyas,
On initial and boundary value problems for fractional order mixed type functional differential inclusions,
{\it Comput. Math. Appl.} 59 (2010) 1253-1265.

\bibitem{dkfa} K. Diethelm, A.D. Freed,
On the solution of nonlinear fractional order differential equations used in the modeling of viscoelasticity, in: F. Keil, W. Mackens, H. Voss, J. Werther (Eds.), Scientific Computing in Chemical Engineering II-Computational Fluid Dynamics, Reaction Engineering and Molecular Properties,
Springer-Verlag, Heidelberg, 1999, 217-224.

\bibitem{jerw} V.J. Ervin, J.P. Roop,
Variational formulation for the stationary fractional advection dispersion equation,
{\it Numer. Meth. Part. Diff. Eqs.} 22 (2006) 58-76.

\bibitem{12} G.J. Fix, J.P. Roop,
Least squares finite-element solution of a fractional order two-point boundary value problem,
{\it Comput. Math. Appl.} 48 (2004) 1017-1033.

\bibitem{3} W.G. Glockle, T.F. Nonnenmacher,
A fractional calulus approach of self-similar protein dynamcs,
{\it Biophys. J.} 68 (1995) 46-53.

\bibitem{hilf} R. Hilfer,
Applications of Fractional Calculus in Physics,
World Scientific, Singapore, 2000.

\bibitem{wj} W. Jiang,
The existence of solutions to boundary value problems of fractional differential equations at resonance,
{\it Nonlinear Anal.} 74 (2011) 1987-1994.

\bibitem{fjy} F. Jiao, Y. Zhou,
Existence results for fractional boundary value problem via critical point theory,
{\it Internat. J. Bifur. Chaos} 22 (2012) Article ID 1250086.

\bibitem{15} A.A. Kilbas, H.M. Srivastava, J.J. Trujillo,
Theory and Applications of Fractional Differential Equations,
Elsevier, Amsterdam, 2006.

\bibitem{kjfx} J.W. Kirchner, X. Feng, C. Neal,
Fractal stream chemistry and its implications for contaminant transport in catchments,
{\it Nature} 403 (2000) 524-526.

\bibitem{lsl} L.S. Leibenson,
General problem of the movement of a compressible fluid in a porous medium,
{\it Izvestiia Akademii Nauk Kirgizsko\u{\i} SSSR} 9 (1983) 7-10.

\bibitem{main} F. Mainardi,
Fractional calculus: some basic problems in continuum and statistical mechanics,
in: A. Carpinteri, F. Mainardi (Eds.), Fractals and Fractional Calculus in Continuum Mechanics,
Springer-Verlag, Wien, 1997, 291-348.

\bibitem{dl1} J. Mawhin, M. Willem,
Critical point theory and Hamiltonian systems,
Applied Mathematical Sciences, Springer, New York, 1989.

\bibitem{dl2} P. Rabinowitz,
Minimax methods in critical point theory with applications to differential equations,
in: CBMS Regional Conference Series in Mathematics, American Mathematical Society, Providence, RI, 1986.

\bibitem{18} S.G. Samko, A.A. kilbas, O.I. Marichev,
Fractional Integrals and Detivatives: Theory and Applications,
Gordon and Breach, New York, 1993.

\bibitem{js} J. Simon,
R\'{e}gularit\'{e} de la solution d'un probl\`{e}me aux limites non lin\'{e}aires,
{\it Ann. Fac. Sci. Tolouse} 3 (1981) 247-274.

\bibitem{zzr} Z. Zhang, R. Yuan,
Infinitely-many solutions for subquadratic fractional Hamiltonian systems with potential changing sign,
{\it Adv. Nonlinear Anal.} 4 (2015) 59-72.

\end{thebibliography}
\end{document}